\documentclass{article}

\newcommand{\oM}{\overline{\mathcal{M}}}
\newcommand{\oP}{\overline{P}}
\newcommand{\oO}{\overline{\Omega}}
\newcommand{\omu}{\overline{\mu}}
\newcommand{\onu}{\overline{\nu}}
\newcommand{\opi}{\overline{\pi}}
\newcommand{\otau}{\overline{\tau}}
\usepackage{amsmath,amsthm,thm-restate,thmtools,amssymb}
\newtheorem{theorem}{Theorem}
\newtheorem{lemma}[theorem]{Lemma}

\newtheorem{example}[theorem]{Example}
\newtheorem{definition}[theorem]{Definition}

\usepackage{graphicx}
\usepackage{nicefrac}
\usepackage{authblk}
\begin{document}

\title{Smaller universes for sampling graphs with fixed degree sequence}
\author[1]{A. Berger}
\author[2]{C. J. Carstens}
\affil[1]{Institute of Computer Science, Martin Luther University Halle-Wittenberg, Halle (Saale), Germany}
\affil[2]{Korteweg-de Vries Institute for Mathematics, University of Amsterdam, Amsterdam, The Netherlands}

\maketitle

\begin{abstract}
{An important problem arising in the study of complex networks, for instance in community detection and motif finding, is the sampling of graphs with fixed degree sequence. The equivalent problem of generating random 0,1 matrices with fixed row and column sums is frequently used as a quantitative tool in ecology. It has however proven very challenging to design sampling algorithms that are both fast and unbiased. 

This article focusses on Markov chain approaches for sampling, where a close-to-random graph is produced by applying a large number $N$ of small changes to a given graph. Examples are the switch chain and Curveball chain, which are both commonly used by practitioners as they are easy to implement and known to sample unbiased when $N$ is large enough. Within theoretical research, much effort has gone into proving bounds on $N$. However, existing theoretical bounds are impractically large for most applications while experiments suggest that much fewer steps are needed to obtain a good sample. 

The contribution of this article is twofold. Firstly it is a step towards better understanding of the discrepancy between experimental observations and theoretically proven bounds. In particular, we argue that while existing Markov chain algorithms run on the set of all labelled graphs with a given degree sequence, node labels are unimportant in practice and are usually ignored in determining experimental bounds. We prove that ignoring node labels corresponds to projecting a Markov chain onto equivalence classes of isomorphic graphs and that the resulting projected Markov chain converges to its stationary distribution at least as fast as the original Markov chain. Often convergence is much faster, as we show in examples, explaining part of the difference between theory and experiments. The speed-up comes from the fact that the projected Markov chain starts from a better (more uniform) initial distribution and runs on a smaller state space. 

Secondly, when node labels \emph{are} important, we show that faster sampling from the space of labelled graphs can be achieved by combining the projected Markov chain with a simple linear-time preprocessing step. As a result of independent interest, we prove that this approach removes the need to introduce `hexagonal moves' in the switch and Curveball chain when sampling labelled simple directed graphs. }
\end{abstract}

\section{Introduction}
The uniform sampling of graphs with fixed degree sequence has attracted a large research effort in network science~\cite{MolloyReed1995,NewmanStrogatzWatts2001,Artzy-Randrup2005,Rao1996}. Samples of random networks are used to determine the significance of properties of real-world networks. For instance, to study the clustering in social networks \cite{Roberts2000}, to understand which subgraphs form the important building blocks of a network~\cite{Milo2002} or to find out if a network is expected to be connected given its degree sequence~\cite{NewmanStrogatzWatts2001}. The equivalent problem of uniformly sampling a 0,1 matrix with fixed row and column sums is considered as one of the most useful `null model approaches' in ecology~\cite{gotelli1996null,Strona2014}. 

Markov chain algorithms, such as the switch chain and Curveball chain, are a popular approach to the above sampling problem \cite{Rao1996,Kannan1999,Strona2014}. Here a graph is randomised by repeatedly applying small degree preserving changes. These algorithms are known to converge to the uniform distribution on the set of \emph{labelled} graphs with a given degree sequence \cite{Artzy-Randrup2005,Rao1996,Verhelst2008,Carstens2015}. Even though Markov chain algorithms are easy to implement and provide a flexible sampling framework, they have a serious drawback: in general it is unknown how many changes need to be applied to obtain a close-to-random graph.  

Two completely separate communities appear to work on the Markov chain approach to this problem. The first one provides us with important theoretical insights to understand the problem in greater depth. It finds graph classes for which the Markov chain can be proven to be efficient~\cite{Cooper07,Greenhill2011,Miklos2013,Erdos2017,Greenhill2018,Amanatidis2018}. Unfortunately the polynomial upper bounds for the running time possess such large exponents, that they can never be used to draw a graph at random in a real world scenario. Furthermore, the graph classes considered (matrices with identical sums for rows or columns) barely occur in the real world.  

The second community uses implementations of these algorithms in practice~\cite{Milo2002,Artzy-Randrup2005,Strona2014,NewmanStrogatzWatts2001}, i.e. they need them for their research to create null models. Applied researchers often stop the Markov chains after a fixed number of steps, using some assessment to judge that this number was large enough such that the sampling was done according to an almost uniform distribution.  

Even though existing theoretical results give impractically large limits, we are optimistic about the speed of this class of algorithms. Several experiments~\cite{Rechner2016,Strona2014,Carstens2017,Ray2015} indicate that both the well known \emph{switch chain}~\cite{Ryser1957,Taylor1981,Kannan1999,Diaconis1995,Rao1996} and the lesser known \emph{Curveball algorithm}~\cite{Verhelst2008,Strona2014,Carstens2015} are quite fast. The only problem is that we can not prove how good they really are, i.e. we have a lack of theory.

In this article we offer a partial explanation for the discrepancy between theoretical and experimental bounds. In most applications, the network statistic of interest only depends on the structure of the network, i.e. it is a \emph{topological} property. In practice, convergence of the statistic of interest is used as an indicator that a Markov chain has converged to its stationary distribution~\cite{Artzy-Randrup2005,Ray2015}. This approach ignores node labels when judging the convergence of the Markov chain. We formally show that ignoring node labels corresponds to projecting a Markov chain onto equivalence classes of isomorphic graphs. We prove that the projected chain converges at least as fast as the original chain and give several examples where convergence is much faster. The speed up is due to sampling from a (often much) smaller state space. In some applications node labels are important, for instance when determining the number of expected edges between certain individuals or communities. We show that faster sampling can be achieved by combining the projected Markov chain with a linear-time preprocessing step. The resulting improved run-time is of clear benefit to practitioners. Furthermore, our contribution opens new pathways for theoretical research, by reducing the size of the state space. As a result of independent interest, we prove that this combination of a projected Markov chain and preprocessing step results in an ergodic Markov chain for all directed graphs, that is, it removes the need for `hexagonal moves'\cite{Rao1996} (Theorem \ref{thm:hex_move}).  

The remainder of this article is organised as follows. We start with a general description of projected Markov chains and prove that the mixing time of these chains is smaller than or equal to the mixing time of the original chain. We then briefly discuss well-known Markov chains used for the sampling of graphs: the switch and Curveball chain. We show that these Markov chains can be projected onto isomorphism classes of graphs and that this results in faster mixing. Furthermore we introduce a preprocessing step which allows us to speed-up the switch and Curveball chain. We give several explicit examples of the method. Finally we discuss limitations and potential extensions to this framework. 

\section{Applications}
The following examples illustrate the main idea behind our speed-up of the switch and Curveball Markov chains. 

\begin{figure}[!htb]
\centering
\includegraphics[width=300px]{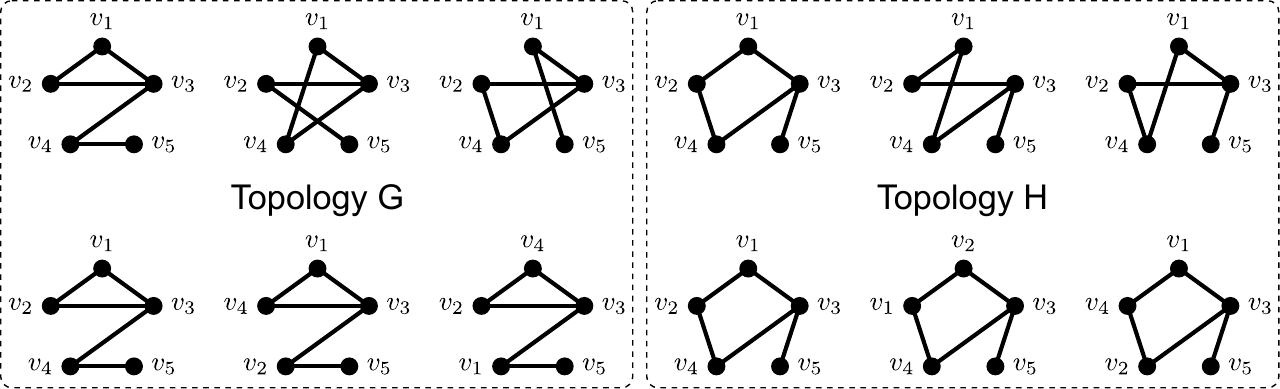}
\caption{\label{fig:example_labels}In the top-row the six simple undirected graphs with degree sequence $(2,2,3,2,1)$, the nodes are fixed in place. There are two sets of three graphs with the same topology. The bottom-row shows the same six graphs, but here the different graphs are illustrated as a relabelling of the nodes.} 
\end{figure}

\begin{example}
Motif finding is a popular tool in network analysis~\cite{Milo2002}. A motif is defined as a small subgraph which appears significantly more frequently in an observed (real-world) graph than in randomly generated graphs. The switch chain is often used to generate such random graphs. 
It samples a graph uniformly at random from the space of all graphs with a given degree sequence. As a small example, Figure \ref{fig:example_labels} illustrates the six different graphs with degree sequence $(2,2,3,2,1)$. Note that the three graphs on the left have the same topology and the three graphs on the right have a second different topology. For motif finding, it is not necessary to generate a sample from all six labelled graphs. We only need to know the probability with which we sample each of these two classes $G$ and $H$, as this allows us to compute the expected number of occurrences of a given subgraph. For instance in this small example, we find a graph with topology $G$ with probability $\nicefrac{1}{2}$ and we find a graph with topology $H$ also with probability $\nicefrac{1}{2}$. Hence, the expected number of triangles equals $0.5$. 
\end{example}

\begin{figure}[!htb]
\centering
\includegraphics[width=300px]{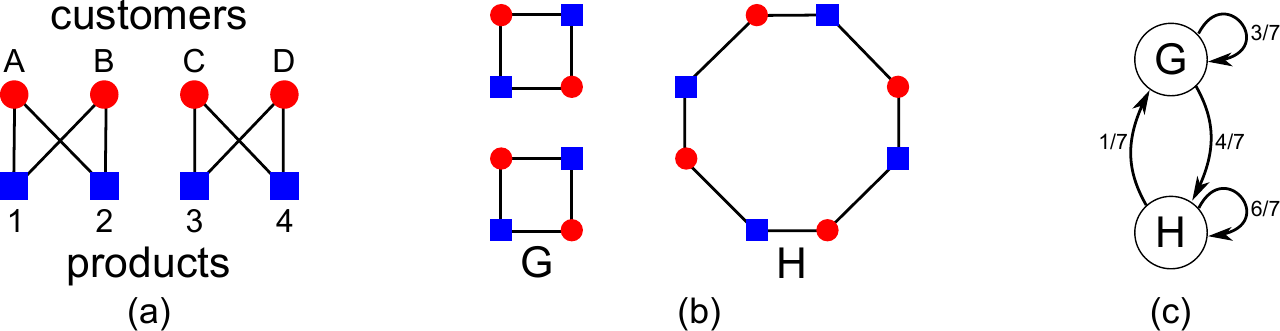}
\caption{\label{fig:example_2222}(a) A toy-example of a customer product network, the customers can be divided into two groups $\{A,B\}$ and $\{C,D\}$ based on the items they have bought. (b) The two unlabelled bipartite graphs with degrees $((2,2,2,2),(2,2,2,2))$. (c) The projected switch chain. It is not hard to see that it has stationary distribution $(\nicefrac{1}{5},\nicefrac{4}{5})$, implying that topology $G$ will be sampled with probability $0.2$.} 
\end{figure}

\begin{example}
\label{exmp:bip_prod_cons}
We are interested in the buying behaviour of customers. We may represent this as a bipartite graph $G$ where the primary nodes represent customers, the secondary nodes represent products and the edges indicate that a customer has bought an item. Say we have observed four customers who each bought two items. In total there are four items and each item has been bought twice. We want to know the probability that the customers can be divided into two groups (see Figure \ref{fig:example_2222}(a)) while fixing the number of items bought per customer and the number of times each item is bought. To do so, we may generate samples of bipartite graphs with the degree sequence $k = ((2,2,2,2),(2,2,2,2))$ using the switch or the Curveball chain, and estimate this probability. Both algorithms sample a graph uniformly at random (provided we run them for long enough) from the set of $90$ distinct \emph{labelled bipartite graphs} with degree sequence $k$. Only $18$ of these states correspond to the situation where we can split the customers into two types. Hence we will find a probability close to $0.2$ provided we take a large enough sample and run the chains for long enough. 

But the property of interest, if the customers can be split into two groups, is a topological property and does not depend on the labelling of the nodes. So in fact we are interested in sampling from a much smaller state space, that of unlabelled bipartite graphs with the given degrees. When removing the node labels, we find that there are only two distinct graphs (as illustrated in Figure \ref{fig:example_2222}(b)). We will later see that we can obtain the probabilities of sampling either of these two topologies by projecting the switch or Curveball chain. After projecting, fewer switches and trades are required to get close to the stationary distribution, largely due to the reduced number of graphs we are sampling from. Figure \ref{fig:example_2222}(c) illustrates this by showing the Markov chain of the projected switch chain. 
\end{example}

These two applications show that it is often unnecessary to sample from the set of labelled graphs. Instead only the topology of the sampled networks is important. The next example shows that even when the labels are important, sampling can be sped up by making use of the projection mechanism. 

\begin{example}
\label{exmp:conn_specific_individuals}
We are studying a social network and want to know what the probability of two specific individuals being connected is given the number of connections of all individuals in the network. In this case, the labels of the nodes, i.e. \emph{who they represent} is important. However, we can still benefit from the speed-up of sampling isomorphic graphs. To see this, we return to the toy-example in Figure \ref{fig:example_labels}. If we generate a sample with the projected switch chain we obtain a given representative of class $G$, say $G_1$, roughly half the time and a representative of $H$, say $H_1$, the rest of the time. To obtain a uniform sample from the six labelled graphs we can use this sample and apply a random node relabelling to all of the sampled graphs. Note that this relabelling has to preserve the degrees of the nodes. Hence in this example we simply choose a random permutation of the node labels $v_1, v_2$ and $v_4$. Now we obtain a random sample from \emph{all} labelled graphs with degree sequence $(2,2,3,2,1)$.  
\end{example}

\section{Projected Markov chains}\label{sec:MC_proj}
We now introduce the framework of projected Markov chains. We show that the projection of a Markov chain has two equivalent interpretations. Firstly we can think of the projection as a Markov chain on equivalence classes. That is the chain $X_0, X_1, \dots, X_t$ becomes $\overline{X_0}, \overline{X_1}, \dots, \overline{X_t}$ and the state space is reduced in size: $\oO := \Omega \backslash \sim$. Secondly, we can interpret the projected chain as running the original chain with an alternative starting distribution: instead of starting in a single state the chain starts from the uniform distribution on states in a single equivalent class. 

In this article we consider discrete time Markov chains $\mathcal{M}=(\Omega, P)$ with finite state space $\Omega$ and transition matrix $P$. It is well-known that such a Markov chain converges to a unique stationary distribution $\pi$ if it is ergodic. If the Markov chain is time-reversible then this distribution satisfies $\pi_i P_{ij} = \pi_j P_{ji}$ for all $X_i,X_j \in \Omega$. All Markov chains discussed in this article are finite, ergodic and time-reversible. 

We will be interested in Markov chains that satisfy condition (\ref{eq:MC_nice_P}) as defined below, and their \emph{projected} or \emph{lumped} chain. We use the following result of Wilson~\cite[Theorem 2.5]{Levin2009}.

\begin{lemma}[Projected Markov chains]\label{thm:projectedChains}
 Let $\mathcal{M}=(\Omega, P)$ be a Markov chain and let $\sim$ be an equivalence relation on $\Omega$ with equivalence classes $[x] \in \bar{\Omega}$. Assume that $P$ satisfies 
 \begin{equation}
P_{x[y]} = P_{x'[y]} \label{eq:MC_nice_P} 
\end{equation}
whenever $x \sim x'$, and where $P_{x[y]}:=\sum_{z\in [y]}P_{xz}$. Then $\bar{\mathcal{M}}=(\bar{\Omega}, \bar{P})$ with $\bar{P}_{[x][y]} := P_{x[y]}$ is a Markov chain. $\bar{\mathcal{M}}$ is called the \emph{projected chain}. 
\end{lemma}

The stationary distribution of the projected Markov chain is proportional to the sizes of the equivalence classes. The following lemma was proved by Grone et al. \cite[Propostion 3]{grone2008interlacing}. 

\begin{lemma} The projected chain $\overline{\mathcal{M}}=(\overline{\Omega}, \overline{P})$ satisfies $\pi_{[x]} P_{[x][y]} = \pi_{[y]} P_{[y][x]}$ where $\pi_{[x]} = \sum_{x\in[x]} \pi(x)$ and hence has stationary distribution $\overline{\pi} = (\pi_{[x_1]}, \dots, \pi_{[x_n]})$. \label{lem:proj_stationary}
\end{lemma} 

The mixing time of a Markov chain quantifies how quickly the chain approaches its stationary distribution. It is defined in terms of the variation distance between distributions.

\begin{definition}
Let $\mu, \nu: \Omega \rightarrow [0,1]$ be probability distributions on $\Omega$. Their variation distance is defined as 
\[d_V(\mu,\nu) = \max_{A \subset \Omega} |\mu(A) - \nu(A)|\]
and can be shown to equal $\nicefrac{1}{2} \sum_{x \in \Omega} |\mu(x) - \nu(x)|$.
\end{definition}

Let $P^t_x$ be the distribution of the Markov chain at time $t$ when started from state $x$, that is when started from the one-point distribution $\mathbf{1}_x$. This distribution $\mathbf{1}_x(y)$ equals $1$ when $x=y$ and $0$ otherwise. When the complete transition matrix $P$ is known, the distribution $P^t_x$ can be computed by $t$ times right multiplying $\mathbf{1}_x$ with $P$, i.e. $P^t_x =  \mathbf{1}_x P^t$. The mixing time of a Markov chain is defined as 

\[\tau(\epsilon) = \max_{x \in \Omega} \min_{T} \{T | d_V(P^t_x, \pi) \leq \epsilon \mbox{ for all } t > T\}.\]

Informally, the mixing time is the maximum number of steps needed to get within distance $\epsilon$ of the stationary distribution regardless of the starting state. We now show that the mixing time $\overline{\tau}(\epsilon)$ of a projected chain is smaller or equal to the mixing time $\tau(\epsilon)$ of the original chain. 

\begin{lemma}\label{thm:ori_comp_proj} Let $\mathcal{M}=(\Omega,P)$ be a finite, ergodic Markov chain with stationary distribution $\pi$ and satisfying (\ref{eq:MC_nice_P}). Then $\tau(\epsilon) \geq \overline{\tau}(\epsilon)$. 

\begin{proof}
We will show that for any $x \in \Omega$, $t \in \mathbb{N}$ we have $d_V(P^t_x, \pi) \geq d_V(\overline{P}^t_{[x]}, \overline{\pi})$ which gives the result. 
Let $f: \Omega \rightarrow \overline{\Omega}$ be the function that maps a state $x$ to its equivalence class $[x]$. Furthermore, for a probability distribution $\mu$ on $\Omega$ let $\mu f^{-1}$ be the probability distribution on $\overline{\Omega}$ given by: 
\[(\mu f^{-1})([x]) := \mu(f^{-1}([x])) = \sum_{z \in [x]} \mu(z).\]
Notice that $\overline{\pi}([x]) = \sum_{x \in [x]} \pi(x) = \pi f^{-1}([x])$, i.e. $\overline{\pi}$ equals $\pi f^{-1}$. Furthermore the `one-point' starting distribution $\overline{P}_{[x]}^0 = \mathbf{1}_{[x]}$ equals $\mathbf{1}_x f^{-1}$ for any $x \in [x]$, that is $\overline{P}_{[x]}^0 = P_x^0 f^{-1}$. We now show that if $\overline{\mu} = \mu f^{-1}$ then also $(\omu \oP) = (\mu P) f^{-1}$. 
We evaluate $(\omu \oP)$ on a class $[y]$ and find
\[\omu\oP([y]) = \sum_{[x] \in \oO} \omu([x]) \oP_{[x][y]} = \sum_{[x] \in \oO} \sum_{z \in [x]} \mu(z) \oP_{[x][y]} = \sum_{z \in \Omega} \mu(z) P_{z[y]}\]
where the first equality comes from writing out the matrix multiplication. The second equality uses $\omu = \mu f^{-1}$ and the last equality uses the definition of $\oP$. Next we obtain 
\[\omu\oP([y]) = \sum_{y \in [y]} \sum_{z \in \Omega} \mu(z) P_{zy} = \sum_{y \in [y]} (\mu P)(y) = (\mu P) f^{-1}([y])\]
by using $P_{z[y]} = \sum_{y \in [y]} P_{zy}$ and again recognizing the matrix multiplication. 
Thus we now know that $\oP_{[x]}^t = P_{x}^t f^{-1}$ for all $t$ and $x$. The proof now follows from \cite[Lemma 7.9]{Levin2009} where it is shown that $d_V(\mu, \nu) \geq d_V(\mu f^{-1}, \nu f^{-1})$ for any $\mu$ and $\nu$. 
\end{proof}
\end{lemma}

We may think of the projected chain $\oM$ as the original chain $\mathcal{M}$ started from the uniform distribution on an equivalence class $[x]$. That is with starting distribution
\[\overline{\mathbf{1}_{x}} = \begin{cases} \frac{1}{|[x]|} & \mbox{when } x \in [x] \\ 0 & \mbox{otherwise.} \end{cases}\]
We will denote by $P^t_{\overline{x}}$ the probability distribution of $\mathcal{M}$ at time $t$ with starting distribution $\overline{\mathbf{1}_{x}}$. We now show that the `mixing time', $\hat{\tau}(\epsilon)$, of $\mathcal{M}$ when started from $\overline{\mathbf{1}_{x}}$ is exactly equal to that of $\oM$, i.e. starting $\mathcal{M}$ from $\overline{\mathbf{1}_x}$ is at least as fast as starting it from $\mathbf{1}_x$ (Lemma 3.3.). To do so we define

\[\hat{\tau}(\epsilon) := \max_{[x] \in \oO} \min_{T} \{T | d_V(P^t_{\overline{x}}, \pi) \leq \epsilon \mbox{ for all } t > T\}.\]

\begin{lemma}\label{thm:proj_equivcl} Let $\mathcal{M}=(\Omega,P)$ be a finite, ergodic Markov chain with stationary distribution $\pi$ the uniform distribution and satisfying (\ref{eq:MC_nice_P}). Then $\hat{\tau}(\epsilon) = \overline{\tau}(\epsilon)$. 

\begin{proof} To prove this statement we will show that $d_V(P^t_{\overline{x}}, \pi) = d_V(\overline{P}^t_{[x]}, \opi)$. \\
Let $\omu: \oO \rightarrow [0,1]$ be a probability distribution on $\oO$, we define a probability distribution $g\omu$ on $\Omega$ by $g\omu(x) := \nicefrac{\omu([x])}{|[x]|}$. Now clearly
\[d_V(g\omu, g\onu) = \frac{1}{2} \sum_{x\in \Omega} \left|\frac{\omu([x])}{|[x]|} - \frac{\onu([x])}{|[x]|}\right| = \frac{1}{2} \sum_{[x] \in \oO} |\omu([x]) - \omu([x])| = d_V(\omu, \onu).\]
Furthermore $P^0_{\overline{x}} = \overline{\mathbf{1}_{x}} = g \mathbf{1}_{[x]} = g P^0_{[x]}$. We next show that if a probability distribution $\mu$ on $\Omega$ can be written as $\mu = g\omu$ then $\mu P = g (\omu\oP)$. We derive 
\[\mu P (y) = \sum_{x \in \Omega} \mu(x) P_{yx} = \sum_{[x] \in \oO} \sum_{x \in [x]} \frac{\omu([x])}{|[x]|} P_{yx} = \sum_{[x] \in \oO} \frac{\omu([x])}{|[x]|} P_{y[x]} = \sum_{[x] \in \oO} \frac{\omu([x])}{|[x]|} \oP_{[y][x]} \]
Now using detailed balance for the the projected chain $\frac{|[x]|}{|\Omega|} \oP_{[x][y]} = \frac{|[y]|}{|\Omega|} \oP_{[y][x]}$ we obtain 
\[\mu P (y) = \sum_{[x] \in \oO} \frac{\omu([x])}{|[y]|} \oP_{[x][y]} = \frac{(\omu\oP)([y])}{|[y]|} = g(\omu\oP)([y]).\]
This implies $P^t_{\overline{x}} = g \oP^t_{[x]}$ for all $t$ and $x$ and thus $d_V(P^t_{\overline{x}}, \pi) = d_V(g \oP^t_{[x]}, g\opi) = d_V(\oP^t_{[x]}, \opi)$.  
\end{proof}
\end{lemma}

In this framework we are able to compare the mixing time of Markov chains \emph{directly} in terms of the variation distance. In the theoretical literature the mixing time of a Markov chain is often bounded using the \emph{spectral gap upper and lower bound}. The spectral gap of an ergodic finite Markov chain $\mathcal{M}=(\Omega, P)$ is defined as follows. Let us denote the left eigenvalues of the transition matrix $P$ by $\lambda_N \leq \dots \leq \lambda_1$. It is a classical result that $-1 < \lambda_N, \lambda_1 =1$ and $\lambda_2 < 1$ (see for example \cite[Lemma 12.1.]{Levin2009}.) Let us denote by $\lambda^* := \max\{|\lambda_i| : 1 < i \leq N\}$ the eigenvalue with largest absolute value smaller than 1. The value $1-\lambda^*$ is often called spectral gap. The main effort for bounding the mixing time of a Markov chain often goes into finding an expression or a bound for the spectral gap. We will later show that our method can sometimes be used to find an explicit expression for the spectral gap of the projected chain (i.e. Example \ref{example:quadratic_to_constant}).

\section{Markov chains for sampling graphs}\label{sec:MC}
There are two commonly used Markov chain algorithms designed for the sampling of graphs; we briefly discuss the switch chain $\mathcal{M}^S$ and the Curveball chain $\mathcal{M}^C$. Both exist in several flavours: for the sampling of bipartite graphs, simple directed graphs, directed graphs and simple undirected graphs. We describe the algorithms in terms of bipartite graphs, in some sense the most general case. The bi-adjacency matrix of a bipartite graph is an $n \times n'$ matrix where $n$ is the number of primary nodes and $n'$ the number of secondary nodes. The $(i,j)$-th entry of this matrix equals one if there is an edge between primary node $p_i$ and secondary node $s_j$, otherwise it equals zero. We will describe the switch chain and Curveball chain as algorithms that randomise the bi-adjacency matrix of a bipartite graph while keeping its row and column sums fixed. Note that this corresponds exactly to sampling a bipartite graph with fixed degrees. 

A Markov chain $\mathcal{M} = (\Omega, P)$ is described by its state space $\Omega$ and its transition probabilities $P$. Given a binary matrix $A$ with row and columns sums $k = ((r_1, \dots, r_n), (c_1, \dots, c_{n'}))$, both $\mathcal{M}^S$ and $\mathcal{M}^C$ have as their state space the set of all binary $n \times n'$ matrices with row and column sums $k$. We denote this state space by $\Omega_k$.

In practice, both Markov chains are started from a specific state $X_0$ (a binary matrix) and each transition corresponds to making a small change to the current state $X_i$. The switch chain applies \emph{switches}: replacing a submatrix 
\begin{equation}
\begin{pmatrix}  1 & 0 \\ 0 & 1 \end{pmatrix} \mbox{ by } \begin{pmatrix}  0 & 1 \\ 1 & 0 \end{pmatrix} \mbox{ or vice versa.} \label{eq:Switch}
\end{equation}
The Curveball algorithm applies \emph{trades}: a trade randomly exchanges `ones' between two selected rows. For instance, the rows 
\begin{equation}
\begin{pmatrix} 1 & 1 & 0 & 1 & 0 & 1 & 0 \\ 0 & 0 & 1 & 1 & 1 & 0 & 0 \end{pmatrix} \mbox{ can be replaced by } \begin{pmatrix} 0 & 1 & 1 & 1 & 1 & 0 & 0 \\ 1 & 0 & 0 & 1 & 0 & 1 & 0 \end{pmatrix} \mbox{.} \label{eq:CB}
\end{equation}
The top row had `tradeable ones', i.e. where the bottom row equals zero, in columns 1,2 and 6 before the trade and the bottom row had tradeable ones in columns 3 and 5. A trade corresponds to randomly selecting three columns from these five available columns. In this example columns 2,3 and 5 were chosen. Notice that also columns 1,2 and 3 could be chosen which corresponds to applying a switch. One trade can apply multiple switches at once. 

There are several versions of the switch chain known in the literature \cite{Kannan1999}, here we use the chain which randomly selects two non-zero matrix entries, and applies a switch if possible, i.e. if the $(2\times 2)$-submatrix formed by the rows and columns corresponding to these entries are as in equation (\ref{eq:Switch}). The Curveball chain that we use in this article proceeds by selecting a pair of rows, $i$ and $j$ at random (with probability $\binom{n}{2}^{-1}$) and applying a trade with probability $\binom{s_i+s_j}{s_i}^{-1}$ where $s_i$ is the number of columns where row $A_i$ equals $1$ and row $A_j$ equals $0$ and $s_j$ is the number of columns where row $A_i$ equals $0$ and row $A_j$ equals $1$, see \cite{Verhelst2008,Strona2014} for more information. 

To define the projected switch and Curveball chain we use the framework discussed in Section \ref{sec:MC_proj}. The equivalence relation that we will use is that of \emph{graph isomorphism}, that is graphs are equivalent when they have the same \emph{topology}. Formally, two (bipartite) graphs $G$ and $H$ are isomorphic if there is a bijective map $\sigma$ between their nodes such that edges are preserved, i.e. such that the edge $(\sigma(u), \sigma(v))$ is present in $H$ if and only if the edge $(u,v)$ is present in $G$. Note that, graph isomorphisms are degree preserving maps. In general it can be hard to decide if two graphs are isomorphic, but creating isomorphic graphs is simple: pick a random labelling of the nodes of a given graph $G$ to obtain an isomorphic graph $H$. 

\begin{figure}[!h]
\centering
\includegraphics[width=300px]{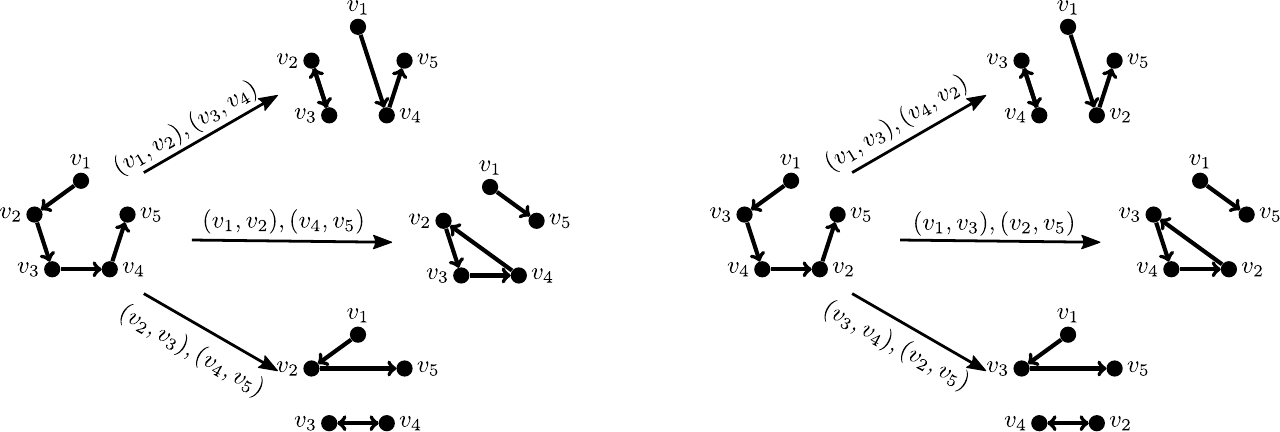}
\caption{\label{fig:iso_nbh}. The switch adjacent graphs of two equivalent graphs $G \sim H$ are identical up to a relabelling of nodes. In fact, the graph isomorphism $\sigma$ between $G$ and $H$ is a graph isomorphism for all graphs that are switch-adjacent to $G$ and $H$. In this case $\sigma$ maps $v_2$ to $v_4$ and $v_4$ to $v_2$ while mapping all other nodes to themselves.} 
\end{figure} 

We say two bipartite graphs $G \sim H \in \Omega$ are equivalent if and only if they are isomorphic as bipartite graphs, i.e. if there exists a graph isomorphism $\sigma = (\sigma_1, \sigma_2)$, where $\sigma_1$ ($\sigma_2$) maps the primary (secondary) nodes of $G$ to the primary (secondary) nodes of $H$. We need to show that the transition matrices $P^S$ and $P^C$ of the switch chain and the Curveball chain respectively are of form (\ref{eq:MC_nice_P}) under this equivalence relation. Intuitively, this holds because applying a specific switch or trade to isomorphic graphs will lead to isomorphic graphs. This implies that the probability of ending up in a given equivalence class is equal for graphs that are in the same equivalence class. In Figure \ref{fig:iso_nbh} we illustrate this for the switch chain of a simple directed graph on five nodes. The following lemma is a formal statement of the above argument in the setting of bipartite graphs, but can easily be generalized to other graph classes. 

\begin{lemma} Let $\mathcal{M}^S=(\Omega,P^S)$, $\mathcal{M}^C=(\Omega,P^C)$ be the switch and Curveball chain for a given bipartite degree sequence. Both $P^S$ and $P^C$ are of the form (\ref{eq:MC_nice_P}). 

\begin{proof}
Let $G=(P,S,E)$ and $G'=(P,S,E')$ be bipartite graphs in $\Omega$ such that $G \sim G'$. That is, there exist degree preserving isomorphisms $\sigma = (\sigma_1, \sigma_2)$, with $\sigma_1: P \rightarrow P$ and $\sigma_2: S \rightarrow S$ such that for any $e=\{p,s\} \in E$ we have $\{\sigma(p),\sigma(s)\} \in E'$. We will write $\sigma(G) := G'$. Let $H$ be switch-adjacent to $G$ with regard to a specific switch: i.e. $H=(V,E \backslash \{\{p_i,s_j\}, \{p_k,s_l\}\} \cup \{\{p_i,s_l\}, \{p_k,s_j\}\}$. Then the graph $\sigma(H) = (V, E'  $ $\backslash \{\{\sigma(p_i),\sigma(s_j)\}, \{\sigma(p_k),\sigma(s_l)\}\}$ $\cup \{\{\sigma(p_i),\sigma(s_l)\}, \{\sigma(p_k),\sigma(s_j)\}\})$ is switch-adjacent to $\sigma(G)$ and furthermore $H \sim \sigma(H)$ by definition. That is, all graphs that are switch-adjacent to $G$ are isomorphic to graphs that are switch adjacent to $\sigma(G)$ under the graph isomorphism $\sigma$ and $P^S_{GH} = P^S_{\sigma(G)\sigma(H)}$. We now write $N_{G}$ for the set of graphs that is switch-adjacent to $G$. 

For any equivalence class $[H]$ and $G \sim G'$ with graph isomorphism $\sigma$ we obtain
\[P^S_{G[H]} = \sum_{K \in [H]}P^S_{GK} = \sum_{K \in N_G \cap [H]}P^S_{GK} \]
\[= \sum_{\sigma(K) \in N_{G'} \cap [H]}P^S_{G'\sigma(K)} = \sum_{\sigma(K) \in [H]}P^S_{G'\sigma(K)} = P^S_{G'[H]}. \]

Since any trade equals a sequence of $k$ switches the result immediately follows for $P^C$. 
\end{proof}
\end{lemma}

In particular, the \emph{projections} $\oM^S$ and $\oM^C$ of the switch chain and the Curveball chain with respect to $\sim$ are well-defined. Lemma \ref{lem:proj_stationary} now tells us that the stationary distribution of these projected chains is proportional to the size of the equivalence classes. Hence, if we generate a sample using the projected chains we obtain each topology with the correct probability: the probability of sampling a graph in $\Omega$ with the given topology. In practice, when we are running experiments where we are only interested in the topology of the sampled networks, we could argue that we already use the projected chain. To illustrate this, we elaborate on Example \ref{exmp:bip_prod_cons} with respect to the switch chain. 

\begin{example}
In Example \ref{exmp:bip_prod_cons} we wanted to know the probability that a bipartite graph with degrees $k = ((2,2,2,2), (2,2,2,2))$ is disconnected. Due to the small size of $\Omega$ and $\oO$ we can explicitly compute $P^S$ and $\oP^S$ (see Figure \ref{fig:example_2222}(c)) and hence determine the mixing time for a given $\epsilon$. For $\epsilon=0.001$ we find $\tau(\epsilon)=28$ and $\otau(\epsilon)=6$. This means that after $28$ switches, the probability of obtaining any specific graph $G$ of the $90$ distinct labelled graphs with degrees $k$ is roughly $\frac{1}{90}$. However, the probability of obtaining a graph with topology $G$ is already roughly $\frac{1}{5}$ after $6$ steps in $\overline{P}^s$.

In general, if we know the mixing time of the switch chain theoretically and run the chain $N$ times for $\tau(\epsilon)$ steps to obtain a sample of size $N$, we could be taking much longer than necessary because the property of interest (and any other topological property) already converges after $\otau(\epsilon)$ steps. 

In fact, for any property of interest (motifs, number of connected components, diameter) we may try to project the chain to an even smaller state space. To see this, let the property be given as a function $f:\Omega \rightarrow \mathbb{R}$. We say two states $s,s'$ are equivalent if and only if $f(s)=f(s')$. Hence, $f$ decomposes $\Omega$ in equivalence classes. If this equivalence relation satisfies 3.1 we obtain a projected chain with smaller mixing time than the original chain. Moreover, since many properties of interest are topological measurements, we know that we can always project the chain down to the isomorphism classes (even if the equivalence relation given by $f$ does not satisfy 3.1). In practice the convergence of the switch chain is often estimated through the convergence of the property of interest. That is, the estimated convergence is the convergence of a projected chain. This can explain part of the difference between theoretically proven bounds and experimentally observed bounds. 
\end{example}

In certain applications, e.g. Example \ref{exmp:conn_specific_individuals}, node labels are important. We now show how we can speed up sampling by adding a preprocessing step to the switch or Curveball chain. We define the preprocessing step in terms of the bi-adjacency matrix. 

\begin{definition}[Preprocessing step]
\label{def:preproc}
Let $A$ be a binary $(n\times n')$-matrix in $\Omega_k$ for fixed row and column sums $k$. Let $R$ be the set of all permutations $\rho:$ $(1, \dots, n)$ $\rightarrow$ $(1, \dots, n)$ such that the row sums $\sum_{j=1}^{n'} A_{\rho(i)j} = r_i$ for all $i$ and let $S$ be the set of all permutations $\sigma:$ $\{1, \dots, n'\}$ $\rightarrow$ $\{1, \dots, n'\}$ such that the column sums $\sum_{i=1}^n A_{i\sigma(j)} = c_j$. We define a preprocessing step which randomly selects a $\rho \in R$ and $\sigma \in S$ to form the matrix $B_{ij} = A_{\rho(i)\sigma(j)}$. This preprocessing step can be implemented by choosing a random order for sets of rows with equal row sum and a random order for sets of columns with equal column sum in $O(n+n')$ \cite{Durstenfeld:1964}.  
\end{definition}

We can think of including the preprocessing step as starting the Markov chain from the uniform distribution on a single equivalence class. That is with starting distribution 
\[\mathbf{1}_{\overline{A}}(B{}) = \begin{cases} \frac{1}{|[A]|} &\mbox{if } B \in [A] \\ 0 & \mbox{otherwise.}\end{cases} \]
Importantly, adding the preprocessing step to the Curveball chain or the switch chain, does not change their convergence to the uniform distribution since the stationary distribution is independent of the starting distribution of the chains. However, it may speed up the convergence. 

For undirected and directed graphs we can similarly introduce a preprocessing step. 

\begin{definition}[Preprocessing step]
\label{def:preproc2}
Let $G$ be a (directed) graph with $n$ nodes and adjacency matrix $A$. Let $P$ be a partition of the nodes with equal degree (equal in- and out-degree in case of directed graphs). Let $R$ be the set of all permutations $\rho:$ $\{1, \dots, n\}$ $\rightarrow$ $\{1, \dots, n\}$ which respect the partition $P$, i.e. nodes are permuted within each partition separately. We define the following preprocessing step: randomly select a permutation $\rho \in R$ and apply it to the nodes of $G$, then return the resulting adjacency matrix. 
\end{definition}

This finishes our discussion of how the framework of projected Markov chains may improve the mixing time in applications. We finish this section by proving that the inclusion of a preprocessing step removes the need for `hexagonal moves' in the switch and Curveball chain for directed graphs as proved by Rao et al~\cite{Rao1996}. The intuition behind this proof is simple: the reason that the state graph of $M^S$ and $M^C$ is disconnected for some directed graphs is that the direction of certain directed triangles can not be reversed~\cite{Berger2010}, which can be achieved with our preprocessing step. 

\begin{theorem}\label{thm:hex_move} The switch chain $M^S$ and Curveball chain $M^C$ have connected state space when including the preprocessing step in Definition \ref{def:preproc2}.

\begin{proof}
It has been shown that the above chains sample uniformly when including a pre-sampling step which assigns orientations to each `induced cycle set' randomly~\cite{Berger2010}. Our suggested preprocessing step does this too, since it permutes nodes with equal in and out degree, and it was shown that three nodes that form an induced cycle set always have the same in- and out-degree~\cite{Lamar2009}. Hence, our preprocessing also re-orients all induced cycle sets. 
\end{proof} 
\end{theorem}

In the next section we will discuss several examples where preprocessing improves the mixing time. Furthermore, we show that for certain families of graphs, the size of the state space $\oO$ is constant whereas the size of $\Omega$ grows quadratically in the number of rows of the matrices. 

\section{Examples: Smaller Universes and Faster Sampling}\label{sec:examples}
In this section we discuss several examples where the state space is reduced significantly by using the framework of projected Markov chains. We start with an example where a quadratically growing state space is reduced to constant size. Furthermore, we show that our method makes it possible to explicitly compute the spectral gap of the projected chain, which would be very complicated for the original chain since its state space is growing in size. The spectral gap is often used to bound the mixing time of Markov chains \cite[Theorem 12.3]{Levin2009}.  
\begin{example}
\label{example:quadratic_to_constant}
Let $r_n$ be a vector of length $n$ where all entries are equal to two, i.e. $r_n=(2,\dots,2)$ and let $c_n$ be the vector $(n-1,n-1,1,1)$. Then $k_n = (r_n, c_n)$ are valid degrees for a bipartite graph with $n$ primary nodes and four secondary nodes $(s_1,s_2,s_3,s_4)$ whenever $n \geq 2$. Let $G$ be the bipartite graph consisting of $K_{n-1,2}$ and $K_{1,2}$, i.e. $s_1$ and $s_2$ are connected to the same $n-1$ primary nodes and $s_3$ and $s_4$ are connected to the same primary node. This graph belongs to the state space $\Omega_{k_n}$. The equivalence class $[G]$ of graph $G$ has size $n$ since there are $n$ choices for the label of the primary node in the disconnected $K_{1,2}$. There is one other equivalence class, the class $[H]$, the graphs in this class are connected, $s_1$ and $s_2$ share $n-2$ neighbours and are each connected to a single additional primary node. One of these two nodes is furthermore connected to $s_3$ and the other to $s_4$. The class $[H]$ has size $2n(n-1)$. Hence the size of the state space $\Omega_{k^n}$ equals $n(2n-1)$ and grows quadratically in terms of $n$. On the other hand, the size of the state space of the projected chain is independent of $n$ and always equal to two.

We are able to explicitly compute the spectral gap for the projected switch chain. The transition probability $P^S_{[G][H]}$ equals the probability of selecting an edge in $K_{n-1,2}$ and an edge in $K_{1,2}$ and hence equals $\binom{m}{2}^{-1} 4 (n-1)$ with $m = 2n$, the total number of edges. The transition probability $P^S_{[H][G]}$ equals $\binom{2n}{2}^{-1} 2$ which can be seen by inspecting a specific graph $H \in [H]$. Let $H$ be the graph where $s_1$ is connected to $p_1, \dots, p_{n-1}$ and $s_2$ is connected to $p_1, \dots, p_{n-2}, p_n$. Furthermore $p_{n-1}$ is connected to $s_3$ and $p_{n}$ is connected to $s_4$. The only switches that will give us a graph in $[G]$ are $(p_{n-1}, s_1)$ with $(p_n, s_4)$ and $(p_{n-1},s_3)$ with $(p_n, s_4)$. 
The eigenvalues of $\oP^S$ can be symbolically computed and equal $\lambda_1=1$, $\lambda_2 = 1 - \nicefrac{2}{n}$. For $n\geq2$, the spectral gap is given by $1-\lambda_2=\nicefrac{2}{n}$, leading to an upper bound of the mixing time of $O(n \log(n/\epsilon))$ \cite{Sinclair1989}.
\end{example}

In the next example we show that for a family of matrices, using our preprocessing step reduces the size of the state space impressively: from growing exponentially in the number of columns to always consisting of a single state. 

\begin{example}
Let $G_l$ be the bipartite graph with two primary nodes with degree $l$ and $2l$ secondary nodes with degree $1$. The size of the state space $\Omega_l$ of graphs with these degrees is $\binom{2l}{l}$ and grows exponentially in $l$. For the Curveball chain, a trade from $G_l$ reaches \emph{all realisations} of the degree sequence with probability $\binom{2l}{l}^{-1}$. Now consider the graphs $G_l$ and $H_l$ where the neighbours of node $p_1$ in $G_l$ are the neighbours of node $p_2$ in $H_l$ and vice versa. To go from state $G_l$ to state $H_l$ with the switch chain, at least $l$ switches are needed. In other words, at least $l$ steps of the switch chain are needed to reach every state with positive probability. This is a clear example where the Curveball chain is the better choice in terms of mixing time. Finally, if we use our preprocessing step and project either chain with respect to the equivalence relation $\sim$, we find that only a single state remains since all states are isomorphic. Hence we only need to apply the preprocessing step and are left with a uniformly sampled labelled graph. 
\end{example}

\section{Discussion and conclusion}\label{sec:discussion}
In this article we introduce a projected version of the switch and Curveball Markov chain where only the topology of the resulting graph is used. In many applications this is the main feature of interest and projecting can significantly reduce the size of the state space and hence improve the mixing time of a Markov chain. We furthermore introduce a preprocessing step that can be used in combination with the projected chain to obtain a random sample from the set of labelled graphs. 

Clearly we can find examples where projecting does not alter the size of the state space, that is, any graph where each node has a unique degree leads to equivalence classes of size one and hence no reduction in the size of the state space. However, such degree sequences have smaller state spaces to begin with, exactly due to the absence of this redundant symmetry. In \cite{Berger2018} an interesting relation is discussed between majorization of degree sequences and the size of the state space. This could turn out to be the key to showing \emph{all} state spaces are small after projection. 

Most theoretical bounds on the mixing time of the switch chain give a bound on the \emph{spectral gap}. The spectral gap of a projected chain is smaller than or equal to the spectral gap of the original chain since it has been proven that the eigenvalues of a projected chain are a subset of the eigenvalues from the original chain \cite[Theorem 12.8.(ii)]{Levin2009}. Focusing the analysis of spectral gap bounds for projected chains would make an interesting area of future research. 

\section*{Acknowledgement}
The authors would like to thank Pieter Kleer for insightful comments and discussions on the content of this work. 


%


\end{document}